\documentclass[12pt,twoside]{article}
\usepackage{amsmath,amsbsy,amsfonts}
\usepackage{color}

\newtheorem{lem}{Lemma}[section]

\newtheorem{theo}{Theorem}[section]

\newtheorem{coro}{Corollary}[section]
\newtheorem{definition}{Definition}[section]

\newtheorem{claim}{Claim}[section]

\let\Section=\section
\def\section{\setcounter{equation}{0}\Section}
\begin{document}
\title{Strauss- and Lions-type results for a class of Orlicz-Sobolev spaces and applications}

\author{Claudianor O. Alves\footnote{C.O. Alves was partially supported by INCT-MAT, PROCAD, CNPq/Brazil
620150/2008-4 and 303080/2009-4}, \\
\noindent Universidade Federal de Campina Grande, \\
\noindent Unidade Acad\^emica de Matem\'atica,\\
\noindent CEP:58429-900, Campina Grande - PB, Brazil\\
\noindent e-mail: coalves@dme.ufcg.edu.br, \\
\vspace{0.5 cm}\\
Giovany M. Figueiredo\thanks{Supported by CNPq/Brazil 300705/2008-5 } \\
\noindent Universidade Federal do Par\'a, \\
\noindent Faculdade de Matem\'atica, \\
\noindent CEP: 66075-110, Bel\'em - Pa, Brazil \\
\noindent e-mail: giovany@ufpa.br \\
\vspace{0.5 cm}\\
and \\
\vspace{0.5 cm}\\
Jefferson A. Santos \\
\noindent Universidade Federal de Campina Grande, \\
\noindent Unidade Acad\^emica de Matem\'atica,\\
\noindent CEP:58429-900, Campina Grande - PB, Brazil\\
\noindent e-mail: jefferson@dme.ufcg.edu.br \\
}

\date{}

\pretolerance10000

\maketitle

\begin{abstract}
The main goal of this work is to prove Strauss- and Lions-type results for Orlicz-Sobolev spaces. After, we use these results to study the existence of solutions for a class of quasilinear problems in $\mathbb{R}^{N}$.

\vspace{0.5 cm}



\end{abstract}

\section{Introduction}

In recent years, a special attention has been given for quasilinear problems of the type
$$
\left\{
\begin{array}{l}
- \mbox{div}(a(|\nabla u|)\nabla u) +  V(x)a(|u|)u=f(u) \ \mbox{in} \ \mathbb{R}^N \\
\mbox{}\\
u \in W^{1}L_{A}(\mathbb{R}^N)  \,\,\, \mbox{with} \,\,\, N \geq 2,
\end{array}
\right.
\eqno{(P)}
$$
where $V,f$ are continuous functions satisfying some technical  conditions and $a:[0,+\infty) \to [0,+\infty)$ is a $C^1$-function.

We cite the papers of Bonanno, Bisci and Radulescu \cite{BBR, BBR2}, Cerny \cite{Cerny}, Cl\'ement, Garcia-Huidobro and Man\'asevich \cite{VGMS}, Donaldson \cite{Donaldson}, Fuchs and Li \cite{Fuchs1},  Fuchs and Osmolovski \cite{Fuchs2}, Fukagai, Ito and Narukawa \cite{Fukagai1,Fukagai2}, Gossez \cite{Gossez}, Le and Schmitt \cite{LK}, Mihailescu and Radulescu \cite{MR1, MR2}, Mihailescu and Repovs \cite{MD}, Mihailescu, Radulescu and Repovs \cite{MRR}, Orlicz \cite{O},  Santos \cite{tesedojefferson} and references therein, where quasilinear problems like $(P)$ have been considered in bounded and unbounded domains of $\mathbb{R}^{N}$. In some those papers, the authors have mentioned that this class of problem arises in a lot of applications, such as, nonlinear elasticity, plasticity and non-Newtonian fluids.

One of the most famous methods to get a solution for $(P)$ is the variational method, where the weak solutions for $(P)$ are precisely the critical points of the energy functional $J: X \to \mathbb{R}$ associated with $(P)$, given by
$$
J(u) =\displaystyle \int_{\mathbb{R}^N}A(|\nabla u|) +
          \displaystyle \int_{\mathbb{R}^N}V(x)A(|u|)
          -\displaystyle \int_{\mathbb{R}^N}F(u),
$$
where $X$ is a convenient subspace of $W^{1}L_{A}(\mathbb{R}^N)$, which depends of the hypotheses on the potential $V$.

In \cite{Fukagai1}, Fukagai, Ito and Narukawa have used the variational method to show the existence of a solution for $(P)$ by assuming that the function $a$ satisfies  the following assumptions: \\

\noindent The function $a(t)t$ is increasing in $(0,+\infty)$, that is,
$$
(a(t)t)'>0 \,\,\, \forall t>0.
\eqno{(a_{1})}
$$

\noindent There exist $l, m \in (1,N)$ such that
$$
l\leq \displaystyle\frac{a(|t|)t^{2}}{A(t)}\leq m \,\,\, \forall t \not= 0,
\eqno{(a_{2})}
$$
where $A(t)=\displaystyle\int^{|t|}_{0}a(s)s \, ds$, $l\leq m < l^{*}$, $l^{*}=\displaystyle\frac{lN}{N-l}$ and $m^{*}=\displaystyle\frac{mN}{N-m}$. \\

Using these hypotheses, the authors showed that $A$ is a N-function satisfying the $\Delta_{2}$ - condition. Moreover, in that paper, it is mentioned some examples of functions $A$, whose function $a(t)$ satisfies the conditions $(a_1)-(a_2)$. The examples are the following
$$
\begin{array}{l}
i) \,\, A(t)=|t|^{p} \,\,\, \mbox{for} \,\,\, 1<p<N. \\
\mbox{}\\
ii) \,\, A(t)=|t|^{p}+|t|^{q} \,\,\, \mbox{for} \,\,\, 1<p<q<N \,\,\, \mbox{and} \,\,\, q \in (p,p^{*}) \,\,\, \mbox{with} \,\,\, \displaystyle p^{*}=\frac{Np}{N-p}.  \\
\mbox{}\\
iii) \,\, A(t)=(1+|t|^2)^{\gamma}-1 \,\,\, \mbox{for} \,\,\, \gamma \in (1, \frac{N}{N-2} ). \\
\mbox{}\\
iv) \,\,  A(t)=|t|^{p}ln(1+|t|) \,\,\, \mbox{for} \,\,\, 1<p_0<p<N-1 \,\,\, \mbox{with} \,\,\, \displaystyle p_0=\frac{-1+\sqrt{1+4N}}{2}.
\end{array}
$$\\

Motivated by \cite{Fukagai1}, more precisely, by hypotheses $(a_1)-(a_2)$ considered on function $a$, the main goal of the present paper is to prove that some results found in Strauss  \cite{Strauss}  and Lions \cite{Lions} also hold in the Orlicz-Sobolev $W^{1}L_A(\mathbb{R}^{N})$  for $A(t)=\int_{0}^{|t|}a(s)s\,ds$, when the above conditions are assumed on $a$. Moreover, results of compactness have been proved for domains in $\mathbb{R}^{N}$, which are invariant by group $O(N)$.

It is well known in the literature, that if the energy  functional is invariant by rotations, sometimes it is possible to find radial solutions for $(P)$. In this case, Strauss-type results can be an interesting tool. Once that we did not find in the literature a Strauss-type result  for Orlicz-Sobolev spaces, the first result of this article goes in this direction and it has the following statement

\begin{theo} \label{ST} ( {\bf A Strauss-type result  for Orlicz-Sobolev spaces } ) \label{strauss}
Assume that $(a_1)-(a_2)$ hold and let $v \in W^{1}L_A(\mathbb{R}^{N})$ be a  radial function. Then
$$
|v(x)| \leq A^{-1}\biggl(\displaystyle\frac{C}{|x|^{N-1}}\displaystyle\int_{\mathbb{R}^{N}}[A(|v|) + A(|\nabla v|)]  \ \biggl)\,\,\, \mbox{a.e in} \,\, \mathbb{R}^{N},
$$
where $A^{-1}$ denotes the inverse function of $A$  restricted to $[0,+\infty)$ and $C$ is a positive constant independent of $v$.
\end{theo}

In the next result, we denote by $ W^{1}L_{A,rad}(\mathbb{R}^{N}) $ the subspace of $ W^{1}L_{A}(\mathbb{R}^{N}) $ consisting of radial functions and by $A_*$ the conjugate function of $A$.

\begin{theo} ({\bf A compactness result for radial functions} ) \label{compacidade} Assume that $(a_1)-(a_2)$ hold and let $B$ be a N-function verifying
$$
\displaystyle\lim_{t\rightarrow 0^+}\displaystyle\frac{B(t)}{A(t)}=0 \eqno{(B_1)}
$$
and
$$
\displaystyle\lim_{t\rightarrow +\infty}\displaystyle\frac{B(t)}{A_{*}(t)}=0.  \eqno{(B_2)}
$$
Then, the embedding $W^{1}L_{A,rad}(\mathbb{R}^{N}) \hookrightarrow L_{B}(\mathbb{R}^{N})$ is compact.
\end{theo}

The above theorem can be applied, when we intend to prove that some functional satisfies, for example, the well known Palais-Smale condition on the space of the radial functions.

In the proof of Theorems \ref{strauss} and \ref{compacidade} the reader is invited to observe that they are true assuming that $A(t)=\int_{0}^{|t|}a(s)s \,ds$ is a $N$-function verifying the $\Delta_2$ condition. Here, we have used conditions $(a_1)-(a_2)$ in view of our applications, see Theorem \ref{gio2} below.

Other important results are  of Lions-type, however we did not find again results of this type for Orlicz-Sobolev spaces. Motivated by this fact, we prove also the following result

 \begin{theo}\label{dolions} ( {\bf A Lions-type result for Orlicz-Sobolev spaces} ) \linebreak Assume that $(a_1)-(a_2)$ hold and let $(u_{n})\subset W^{1}L_{A}(\mathbb{R}^N)$ be a bounded sequence  such that there exists $R>0$ satisfying
$$
\displaystyle \lim_{n\rightarrow +\infty}\displaystyle \sup_{y\in{\mathbb{R}^N}}\int_{B_{R}(y)}A(|u_{n}|)=0.
$$
Then, for any N-function $B$ verifying $\Delta_{2}$-condition with
$$
\displaystyle\lim_{t\rightarrow 0}\displaystyle\frac{B(t)}{A(t)}=0 \eqno{(B_1)}
$$
and
$$
\displaystyle\lim_{|t|\rightarrow +\infty}\displaystyle\frac{B(t)}{A_{*}(t)}=0,  \eqno{(B_2)}
$$
we have
$$
u_{n}\rightarrow 0 \ \ \mbox{in} \ \ L_{B}(\mathbb{R}^{N}).
$$
\end{theo}

Theorem \ref{dolions} is interesting because it can be used to prove the existence of critical points for the energy functional $J$, when the potential $V$ is \linebreak  $\mathbb{Z}^{N}$- periodic.

Our next result can also be used to show compactness results for the space $W^{1}_{0}L_{A}(\Omega)$, when $\Omega \subset \mathbb{R}^{N}$ is invariant with respect to  action of a subgroup of $O(N)$.  Before to state it, we need to fix some definitions and notations. To this end, we follow the spirit of Willem's book \cite{Willem}.

\begin{definition} \label{DEF1} Let $G$ be a subgroup of $O(N), y \in \mathbb{R}^{N}$ and $r>0$. We define,
$$
m(y,r,G)=\sup \{n \in \mathbb{N}: \exists g_1,...,g_n \in G: j \not= k \Rightarrow B_r(g_j y) \cap B_r(g_k y)= \emptyset \}.
$$
\end{definition}

An open set $\Omega \subset \mathbb{R}^{N}$ is said invariant when $g \Omega =\Omega$ for all $g \in G$. An invariant subset $\Omega \subset \mathbb{R}^{N}$ is compatible with $G$ if, for some $r>0$,
$$
\lim_{
\begin{array}{l}
|y| \to +\infty \\
dist(y, \Omega) \leq r
\end{array}}
m(y,r,G)=+\infty .
$$

\begin{definition} \label{DEF2} Let $G$ be a subgroup of $O(N)$ and let $\Omega \subset \mathbb{R}^{N}$ be an invariant set. The action of $G$ on  $W^{1}_{0}L_{A}(\Omega)$ is defined by
$$
gu(x)=u(g^{-1}x) \,\,\, \forall x \in \mathbb{R}^{N}.
$$
The subspace of invariant functions is defined by
$$
W^{1}_{0,G}L_{A}(\Omega)=\{ u \in  W^{1}_{0}L_{A}(\Omega): gu=u, \,\,\, \forall g \in G\}.
$$
\end{definition}

\begin{theo} \label{T5} ( {\bf A Compactness result involving the group $O(N)$.} )  If $\Omega$ is compatible with $G$ and $(a_1)-(a_2)$ hold, the embedding
$$
W^{1}_{0,G}L_{A}(\Omega) \hookrightarrow L_{B}(\Omega),
$$
is compact, for any N-function $B$ verifying $\Delta_{2}$-condition with
$$
\displaystyle\lim_{t\rightarrow 0}\displaystyle\frac{B(t)}{A(t)}=0 \eqno{(B_1)}
$$
and
$$
\displaystyle\lim_{|t|\rightarrow +\infty}\displaystyle\frac{B(t)}{A_{*}(t)}=0. \eqno{(B_2)}
$$

\end{theo}

As an immediate consequence of the last result, we have the following corollary

\begin{coro} \label{C1} Let $N_j \geq 2, j=1,...,k, \sum_{j=1}^{k}N_j=N$ and
$$
G=O(N_1) \times O(N_2) \times .... \times O(N_k).
$$
Then, the compact embeddings of Theorem \ref{T5} occur with $\Omega = \mathbb{R}^{N}$.

\end{coro}

Related to the Theorems \ref{compacidade}, \ref{dolions}, \ref{T5} and Corollary \ref{C1}, we would like to cite the paper due to Fan, Zhao and Zhao \cite{FZZ}, where results like above has been established for the space $W^{1,p(x)}(\mathbb{R}^{N})$.

Motivated by the above results, we study the existence of solutions for some classes of quasilinear problems assuming that  $V:\mathbb{R}^{N}\rightarrow \mathbb{R}$ is a continuous function verifying
$$
0<V_0 =\inf_{x \in \mathbb{R}^{N}}V(x) \eqno{(V_1)}
$$
and $f:\mathbb{R}\rightarrow \mathbb{R}$ is a $C^1$-function satisfying the  properties:
$$
\lim_{|t|\rightarrow 0}\displaystyle\frac{f(t)}{a(|t|)|t|}=0
\eqno{(f_{1})}
$$
and
$$
\lim_{|t|\rightarrow +\infty}\displaystyle\frac{f(t)}{a_{*}(|t|)|t|}=0,
\eqno{(f_{2})}
$$
where $a_{*}(t)t$ is such that the Sobolev conjugate function $A_{*}$ of $A$ ( see Section 2) is its primitive, that is,
$A_*(t)= \displaystyle \int_{0}^{|t|}a_{*}(s)s \, ds$.

There exists $\theta > m$ such that
$$
0<\theta F(t)=\int^{t}_{0}f(s)ds \leq tf(t) \,\,\, \mbox{for all}
\,\,\, t \in \mathbb{R} \setminus \{0\}.
\eqno{(f_{3})}
$$

Our main result concerning the existence of a solution for problem $(P)$ is the following:

\begin{theo}\label{gio2}
Suppose that $(f_{1})-(f_{3})$,  $(a_{1})-(a_{2})$ and $(V_1)$ hold. Moreover, assume that one of the following conditions hold: \\

\noindent $i)$ \,  $V$ is a radial function, that is,
$$
V(x)=V(|x|), \,\,\,\, \forall x \in \mathbb{R}^{N},
$$
\noindent or \\
\noindent $ii)$ \, $V$ is a $\mathbb{Z}^{N}$ periodic function, that is,
$$
V(x+y)=V(x), \,\,\, \forall x \in \mathbb{R}^{N} \,\,\, \mbox{and} \,\,\, \forall y \in \mathbb{Z}^{N}.
$$
Then, problem $(P)$ has a nontrivial solution.
\end{theo}

\vspace{0.5 cm}

The plan of this paper is as follows. In Section 2, we review some proprieties of Orlicz and Orlicz-Sobolev spaces. In Section 3, we prove Theorems \ref{strauss}, \ref{compacidade}, \ref{dolions} and \ref{T5}.  In Section 4, we given a proof of Theorem \ref{gio2}.

\section{ A brief review about $N$-function and Orlicz-Sobolev spaces}

\mbox{}\,\,\,\,\,\,\,\,\, In this section, we recall some properties of Orlicz and Orlicz-Sobolev spaces. The reader can find more properties of these spaces in the books of Adams and Fournier \cite{adams}, Adams and Hedberg \cite{adams2}, Donaldson and Trundiger \cite{Donaldson2}, Fuchs and Osmolovski \cite{Fuchs2}, Musielak \cite{M} and O'Neill \cite{Oneill}.

First of all, we recall that a continuous function $\Phi:\mathbb{R} \to [0,+\infty)$ is a $N$-function if: \\

\noindent $i)$ \, $\Phi$  is convex. \\

\noindent $ii)$ \, $\Phi(t)=0 \Leftrightarrow t=0$. \\

\noindent $iii)$ \, $ \displaystyle \frac{\Phi(t)}{t} \stackrel{t \to 0}{\longrightarrow} 0 \,\,\, \mbox{and} \,\,\, \frac{\Phi(t)}{t} \stackrel{t \to +\infty }{\longrightarrow}+\infty $. \\

\noindent $iv)$ \, $\Phi$ is even.

\vspace{0.5 cm}
In what follows, we say that a $N$-function $\Phi$ verifies the $\Delta_{2}$-condition if
$$
\Phi(2t)\leq K \Phi(t) \,\,\, \forall \,\, t\geq 0,
$$
for some constant $K>0$. This condition can be rewritten in the following way: For each $s >0$, there exists $M_s>0$ such that
$$
\Phi(st) \leq M_s \Phi(t) \,\,\, \forall t \geq 0.    \eqno{(\Delta_2)}
$$

Fixed an open set $\Omega \subset \mathbb{R}^{N}$ and a N-function $\Phi$, the Orlicz space $L_{\Phi}(\Omega)$ is defined. When $\Phi$ satisfies $\Delta_{2}$-condition, the space $L_{\Phi}(\Omega)$ is the vectorial space  of the measurable functions $u: \Omega \to \mathbb{R}$ such that
$$
\displaystyle\int_{\Omega}\Phi(|u|) < \infty.
$$
The space $L_{\Phi}(\Omega)$ endowed with Luxemburg norm, that is, with the norm given by
$$
\|u\|_{\Phi}= \inf \biggl\{\alpha >0: \int_{\Omega}\Phi\Big(\frac{|u|}{\alpha}\Big)\leq 1\biggl\},
$$
is a Banach space. The complement function of $\Phi$, denoted by $\widetilde{\Phi}(s)$, is given by the Legendre transformation, that is
$$
\widetilde{\Phi}(s)=\displaystyle\max_{t \geq 0}\{st -\Phi(t)\} \ \ \mbox{for} \ \ s \geq 0.
$$
The functions $\Phi$ and $\widetilde{\Phi}$ are complementary each other. Moreover, we have the Young's inequality given by
\begin{equation} \label{D2}
st \leq \Phi(t) + \widetilde{\Phi}(s) \,\,\,\,\,\, \forall t,s \geq 0.
\end{equation}
Using the above inequality, it is possible to prove a H\"{o}lder type inequality, that is,
\begin{equation} \label{D3}
\biggl|\displaystyle\int_{\Omega}u v\biggl|\leq 2\|u\|_{\Phi}\|v\|_{\widetilde{\Phi}}, \,\,\, \forall \,\, u \in L_{\Phi}(\Omega) \,\,\, \mbox{and} \,\,\,  v \in L_{\widetilde{\Phi}}(\Omega).
\end{equation}
Another important function related to function $\Phi$, it is the Sobolev conjugate function $\Phi_{*}$ of $\Phi$ defined by
$$
\Phi^{-1}_{*}(t)=\displaystyle\int^{t}_{0}\displaystyle\frac{\Phi^{-1}(s)}{s^{(N+1)/N}} ds \,\,\, \mbox{for} \,\,\, t >0,
$$
when
$$
\displaystyle\int^{+\infty}_{1}\displaystyle\frac{\Phi^{-1}(s)}{s^{(N+1)/N}} ds=+\infty.
$$

If $\Phi(t)=|t|^{p}$ for $1< p<N $, we have $\Phi_{*}(t)={p^{*}}^{p^{*}}|t|^{p^{*}}$, where $p^{*}=\frac{pN}{N-p}$.

The next lemma will be used in the proof of some results and its proof can be found in \cite{Gossez}

\begin{lem}\label{brezislieb} Let $\Omega\subset \mathbb{R}^N$ be an open set and $\Phi:\mathbb{R}\rightarrow [0,\infty)$
be a N-function satisfying the $\Delta_2-$condition. If also the 
complementary function $\widetilde{\Phi}$ satisfies the
$\Delta_2-$condition and $(f_n)$ is a bounded sequence in $L_\Phi(\Omega)$ satisfying 
$$
f_n(x)\rightarrow f(x)\  \text{a.e. in } \, \Omega,
$$
then
$$
f_n\rightharpoonup f \ \text{in }L_\Phi(\Omega),
$$
that is,
$$
\int_{\Omega}f_n v \, dx \to \int_{\Omega}f v \, dx \,\,\, \forall v \in L_{\tilde{\Phi}}(\Omega).
$$

\end{lem}

Hereafter, we denote by $W^{1}_{0}L_{\Phi}(\Omega)$ the Orlicz-Sobolev space obtained by the completion of $C^{\infty}_{0}(\Omega)$ with the norm
$$
\|u\|=\|\nabla u\|_{\Phi}+ \|u\|_{\Phi}.
$$
When $\Omega = \mathbb{R}^{N}$, we use the symbol $W^{1}L_{\Phi}(\mathbb{R}^N)$ to denote the space $W^{1}_{0}L_{\Phi}(\mathbb{R}^{N})$.

An important property that we must detach is: If $\Phi$ and $\widetilde{\Phi} $ satisfy $\Delta_2$-condition, the spaces $L_{\Phi}(\Omega)$ and $W^{1}L_{\Phi}(\mathbb{R}^N)$ are reflexive and separable. Moreover, the $\Delta_2$-condition also implies that
\begin{equation} \label{CV0}
u_n \to u \,\,\, \mbox{in} \,\,\, L_{\Phi}(\Omega) \Longleftrightarrow \int_{\Omega}\Phi(|u_n-u|) \to 0
\end{equation}
and
\begin{equation} \label{CV1}
u_n \to u \,\,\, \mbox{in} \,\,\, W^{1}L_{\Phi}(\Omega) \Longleftrightarrow \int_{\Omega}\Phi(|u_n-u|) \to 0 \,\,\, \mbox{and} \,\,\, \int_{\Omega}\Phi(|\nabla u_n- \nabla u|) \to 0.
\end{equation}

In the literature, we find some important embeddings involving the Orlicz-Sobolev spaces, for example, it is possible to prove that embedding
$$
W^{1}L_{\Phi}(\mathbb{R}^N) \hookrightarrow L_{B}(\mathbb{R}^{N})
$$
is continuous,  if $B$ is a $N$-function satisfying 
$$
\limsup_{t \to 0}\frac{B(t)}{\Phi(t)} < +\infty \,\,\, \mbox{and} \,\,\,  \limsup_{|t| \to +\infty}\frac{B(t)}{\Phi_{*}(t)} < +\infty.
$$
When the space $\mathbb{R}^{N}$ is replaced by a bounded domain $D$ and the limits below hold
\begin{equation} \label{M1}
\limsup_{t \to 0}\frac{B(t)}{\Phi(t)}< +\infty \,\,\, \mbox{and} \,\,\, \limsup_{|t| \to +\infty}\frac{B(t)}{\Phi_{*}(t)}=0,
\end{equation}
the embedding
\begin{equation} \label{M2}
W^{1}L_{\Phi}(D) \hookrightarrow L_{B}(D)
\end{equation}
is compact.

\vspace{0.5 cm}

The next four lemmas involve the functions $A, \widetilde{A}$ and $A_{*}$ and theirs proofs can be found in \cite{Fukagai1}. Hereafter, $A$ is the $N$-function given in the introduction and $\widetilde{A},A_{*}$ are the complement and conjugate functions of $A$ respectively.

\begin{lem} \label{F0} The functions $A$ and $\widetilde{A}$ satisfy the inequality
\begin{equation} \label{D1}
\widetilde{A}(a(|t|)t) \leq A(2t) \,\,\, \forall t \geq 0.
\end{equation}
\end{lem}

\begin{lem} \label{F1} Assume that $(a_1)-(a_2)$ hold and let $\xi_{0}(t)=\min\{t^{l},t^{m}\}$,\linebreak $ \xi_{1}(t)=\max\{t^{l},t^{m}\},$ for all $t\geq 0$. Then,
$$
\xi_{0}(\rho)A(t) \leq A(\rho t) \leq  \xi_{1}(\rho)A(t) \;\;\; \mbox{for} \;\; \rho, t \geq 0
$$
and
$$
\xi_{0}(\|u\|_{A}) \leq \int_{\mathbb{R}^{N}}A(|u|) \leq \xi_{1}(\|u\|_{A})  \;\;\; \mbox{for} \;\; u \in L_{A}(\mathbb{R}^{N}).
$$
\end{lem}

\begin{lem} \label{F3} The function $A_*$ satisfies the following inequality
$$
l^{*} \leq \frac{a_*(|t|)t^{2}}{A_{*}(t)} \leq m^{*} \,\,\, \mbox{for} \,\,\, t \not= 0.
$$
\end{lem}
As an immediate consequence of the Lemma \ref{F3}, we have the following result

\begin{lem} \label{F2} Assume that $(a_1)-(a_2)$ hold and let  $\xi_{2}(t)=\min\{t^{l^{*}},t^{m^{*}}\},$ $\xi_{3}(t)=\max\{t^{l^{*}},t^{m^{*}}\}$, for all $t\geq 0$. Then,
$$
\xi_{2}(\rho)A_*(t) \leq A_*(\rho t) \leq \xi_{3}(\rho)A_*(t) \;\;\; \mbox{for} \;\; \rho, t \geq 0
$$
and
$$
\xi_{2}(\|u\|_{A_*}) \leq \int_{\mathbb{R}^{N}}A_*(|u|) \leq \xi_{3}(\|u\|_{A_*})  \;\;\; \mbox{for} \;\; u \in L_{A_*}(\mathbb{R}^{N}).
$$
\end{lem}

\section{Strauss- and Lions- type results for Orlicz-Sobolev spaces}

After the above brief review, we are able to prove our main results involving the Orlicz-Sobolev spaces.

\vspace{0.5 cm}

\noindent \textbf{Proof of Strauss' Theorem ( Theorem \ref{strauss} ).} First of all, we will establish the  result for functions in $C^{\infty}_0(\mathbb{R}^{N})$. After, by density, we establish the result for all radial functions in  $W^{1}L_A(\mathbb{R}^{N})$.

Consider $v \in C_{0}^{\infty}(\mathbb{R}^N)$, $|x|=r$ and $ w(r)=v(x)$. Note that
$$
A(w(b))-A(w(r))= \displaystyle\int^{b}_{r}\biggl(\displaystyle\frac{d}{ds}A(w)\biggl) \ ds \,\,\, \forall \, b>r>0.
$$
Since $w \in C^{\infty}_{0}([0,\infty))$, for $b$ large enough,
$$
A(w(r)) = - \displaystyle\int^{\infty}_{r} a(|w|)ww' \ ds\leq \displaystyle\int^{\infty}_{r} a(|w|)|w||w'| \ ds.
$$
Combining (\ref{D2}) with (\ref{D1})
$$
a(|w|)|w||w'| \leq \widetilde{A}(a(|w|)|w|) +A(|w'|) \leq A(2 |w|) + A(|w'|),
$$
then by $\Delta_{2}$-condition,
$$
a(|w|)|w||w'| \leq K A(|w|) +A(|w'|).
$$
Therefore,
$$
A(w(r)) \leq (K+1) \displaystyle\int^{\infty}_{r}[A(|w(s)|)+ A(|w'(s)|)] ds,
$$
and we can conclude that
$$
A(w(r)) \leq \frac{(K+1)}{r^{N-1}}\int^{\infty}_{r}[A(|w(s)|)+ A(|w'(s)|)]s^{N-1} ds.
$$
From this, there is $C>0$ such that
$$
A(v(x)) \leq \displaystyle\frac{C}{|x|^{N-1}}\displaystyle\int_{\mathbb{R}^{N}}[A(|v|)+ A(|\nabla v|)].
$$
Since $A$ is an even function, $A(v(x))=A(|v(x)|)$ for all $x \in \mathbb{R}^{N}$, and so,
$$
A(|v(x)|) \leq \displaystyle\frac{C}{|x|^{N-1}}\displaystyle\int_{\mathbb{R}^{N}}[A(|v|)+ A(|\nabla v|)].
$$
From this,
$$
|v(x)| \leq A^{-1}\biggl(\displaystyle\frac{C}{|x|^{N-1}}\displaystyle\int_{\mathbb{R}^{N}}[A(|v|)+ A(|\nabla v|)]\biggl) \,\,\, \forall x \in \mathbb{R}^{N} \setminus \{0\},
$$
where $A^{-1}$ denotes the inverse function of $A$ restricted to $[0,+\infty)$. Now, the lemma follows from the density of $C^{\infty}_{0}(\mathbb{R}^{N})$ in $W^{1}L_{A}(\mathbb{R}^{N})$.
$\rule{2mm}{2mm}$

\vspace{0.5 cm}

Now, we are able to prove the  compactness result involving $ W^{1}L_{A,rad}(\mathbb{R}^{N}) $.

\vspace{0.5 cm}

\noindent \textbf{Proof of the compactness Theorem ( Theorem \ref{compacidade} ).} Let $\{u_n\} \subset W^{1}L_{A,rad}(\mathbb{R}^{N})$ be a sequence verifying
$$
u_n \rightharpoonup 0 \;\;\; \mbox{in} \;\;\; W^{1}L_{A,rad}(\mathbb{R}^{N}).
$$
Without loss of generality, we can assume that $u_n \geq 0$ for all $n \in \mathbb{N}$. From $(B_1)-(B_2)$, for each $\epsilon>0$ and $q >1$, there is $C>0$ such that
\begin{equation}\label{B1}
B(t)\leq \epsilon \left(A(t)+A_*(t)\right)+ C |t|^q \;\; \forall t\geq0.
\end{equation}
Using Theorem \ref{strauss}, Lemma \ref{F1} and the boundedness of $\{u_n\}$  in $W^{1}L_{A}(\mathbb{R}^{N})$, for each $R>0$, there is $C>0$ such that
$$
|u_n(x)|^{q} \leq C \left(\frac{1}{|x|^{\frac{(N-1)}{m}q}} + \frac{1}{|x|^{\frac{(N-1)}{l}q}} \right)\;\;\; \mbox{in} \,\,\, [|x| \geq R]  \,\,\, \mbox{and} \,\,\, \forall n \in \mathbb{N}.
$$
Choosing $q$ large enough,
$$
g(x)= C \left(\frac{1}{|x|^{\frac{(N-1)}{m}q}} + \frac{1}{|x|^{\frac{(N-1)}{l}q}}\right)  \in L^{1}([|x|>\delta]) \;\;\; \forall \delta >0.
$$
The last inequality combined with Lebesgue's Theorem implies that
$$
\int_{[|x| \geq R]}|u_n(x)|^{q} \to 0 \;\;\; \mbox{as} \;\;\; n \to \infty.
$$
This limit together with (\ref{B1}) leads to
\begin{equation} \label{B2}
\int_{[|x|\geq R]}B(u_n) \to 0 \;\;\; \mbox{as} \;\; n \to \infty.
\end{equation}
Observing that $(B_1)-(B_2)$ imply that  $(\ref{M1})-(\ref{M2})$ hold, one has that 
$$
W^{1}L_{A}([|x|<R]) \hookrightarrow L_{B}([|x|<R])
$$
is a compact embedding. Hence,
\begin{equation} \label{B3}
\int_{[|x| <R]}B(u_n) \to 0 \;\;\; \mbox{as} \;\; n \to \infty.
\end{equation}
From (\ref{B2}) and (\ref{B3}),
$$
\int_{\mathbb{R}^{N}}B(u_n) \to 0 \;\;\; \mbox{as} \;\; n \to \infty,
$$
and the proof of the theorem is complete. $\rule{2mm}{2mm}$

\vspace{0.5 cm}

\noindent \textbf{Proof of the Lions' Theorem ( Theorem \ref{dolions} ).} First of all, we observe that
$$
\int_{\mathbb{R}^N}B(|u_n|) =\int_{[\mid u_n\mid>k]}B(|u_n|)+\int_{[\mid u_n\mid\leq k]}B(|u_n|).
$$
From $(B_2)$, given $\epsilon>0$, there is $k>0$ such that
$$
B(t)=B(|t|)\leq \epsilon A_*(| t|), \ \text{if } |t| >k,
$$
which yields
$$
\int_{[|u_n|>k]}B(|u_n|)\leq \epsilon\int_{[|u_n|>k]}A_*(|u_n|) \leq \epsilon C,
$$
and so,
$$
\displaystyle \limsup_{n\rightarrow+\infty} \int_{\mathbb{R}^N}B(|u_n|)\leq \epsilon C+\limsup_{n\rightarrow+\infty}\int_{[\mid u_n\mid\leq k]}B(|u_n|).
$$
\begin{claim} \label{Claim1} $\displaystyle \limsup_{n\rightarrow+\infty}\int_{[\mid u_n\mid\leq k]}B(|u_n|)=0$
\end{claim}
Using this Claim,
$$
\displaystyle \limsup_{n\rightarrow+\infty} \int_{\mathbb{R}^N}B(|u_n|)\leq \epsilon C,
$$
from where it follows that
$$
\displaystyle \limsup_{n\rightarrow+\infty} \int_{\mathbb{R}^N}B(|u_n|)=0,
$$
then
$$
u_n\buildrel  n \to +\infty \over
\longrightarrow0\ \text{em }L_B(\mathbb{R}^N).
$$
Now, we will prove the Claim \ref{Claim1}. Setting the function
$$
v_n(x)=\chi_{[\mid u_n\mid\leq k]}(x)u_n(x),
$$
it is sufficient to show that
\begin{equation}\label{Eq.1}
\displaystyle \limsup_{n\rightarrow+\infty} \int_{\mathbb{R}^N}B(|v_n|)=0.
\end{equation}
From  $(\Delta_2)$, there is $M_k>0$ such that
$$
A\left(\left|\frac{v_n}{k}\right|\right) \leq M_k A(|v_n|), \,\,\, \forall n \in \mathbb{N}.
$$
This combined with Lemma \ref{F1} asserts
$$
\int_{B_R(y)}A(|v_n|) \geq \frac{1}{M_k}\int_{B_R(y)}A\left(\left|\frac{v_n}{k}\right|\right)\geq C \int_{B_R(y)} \left| \frac{v_n}{k}\right|^{m},
$$
and so,
$$
\displaystyle \lim_{n\rightarrow+\infty}\sup_{y \in \mathbb{R}^{N}} \int_{B_R(y)}\left| \frac{v_n}{k}\right|^{m}=0.
$$
Fixing
$$
w_n=\frac{v_n}{k}, (| w_n|_\infty \leq 1),
$$
we get
\begin{equation} \label{WN}
\displaystyle \lim_{n\rightarrow+\infty}\sup_{y \in \mathbb{R}^{N}} \int_{B_R(y)}|w_n|^m=0.
\end{equation}
Using again $(\Delta_2)$, there is  $\widehat{M}_k>0$ such that
$$
\int_{\mathbb{R}^N}B(|v_n|)=\int_{\mathbb{R}^N}B\left(k\frac{|v_n|}{k}\right)\leq \widehat{M}_k \int_{\mathbb{R}^N}B(|w_n|).
$$
Consequently, the limit (\ref{Eq.1}) follows if
$$
\displaystyle \limsup_{n\rightarrow+\infty} \int_{\mathbb{R}^N}B(|w_n|)=0.
$$

\begin{claim} \label{Claim2}
For all $\alpha>1$ and $n \in \mathbb{N}$,  $A(|w_n|^\alpha)\in W^{1,1}(\mathbb{R}^N)$.
\end{claim}
Indeed, since $|w_n|_\infty\leq 1$ and  $w_n\in W^{1}L_A(\mathbb{R}^N)$,
\begin{equation} \label{WW1}
\int_{\mathbb{R}^N}A(|w_n|^\alpha) \leq \int_{\mathbb{R}^N}A(|w_n|) < +\infty  \,\,\, \mbox{and} \,\,\, \int_{\mathbb{R}^N}A(|\nabla w_n|) < +\infty .
\end{equation}
Moreover,
\begin{eqnarray*}
\int_{\mathbb{R}^N}\mid\nabla (A(|w_n|^\alpha))\mid&\leq& \alpha\int_{\mathbb{R}^N}a(| w_n|^\alpha)| w_n|^{\alpha}|w_n|^{\alpha-1}|\nabla w_n|\\
&\leq& \alpha\int_{\mathbb{R}^N}a(| w_n|^\alpha)| w_n|^{\alpha}|\nabla w_n|.
\end{eqnarray*}
Since by (\ref{D2}) and (\ref{D1}),
$$
a(|w_n|^\alpha)|w_n|^{\alpha}|\nabla w_n| \leq \widetilde{A}(a(|w_n|^{\alpha})|w_n|^{\alpha}) +A(|\nabla w_n|) \leq A(2 |w_n|^{\alpha}) + A(|\nabla w_n|),
$$
the $\Delta_{2}$-condition yields,
$$
a(|w_n|^{\alpha})|w_n|^{\alpha}|\nabla w_n| \leq K A(|w_n|^{\alpha}) +A(|\nabla w_n|),
$$
therefore,  (\ref{WW1}) gives
$$
\int_{\mathbb{R}^N}\mid\nabla (A(| w_n|^\alpha))\mid < +\infty.
$$
By Sobolev embedding,
$$
W^{1,1}(B_R(y))\hookrightarrow L^{\frac{N}{N-1}}(B_R(y)).
$$
Therefore, there exists $C>0$ such that
$$
\left(\int_{B_R(y)}A(| w_n|^\alpha)^{\frac{N}{N-1}}\right)^{\frac{N-1}{N}}\leq C \int_{B_R(y)}\left(\left| \nabla A(| w_n|^\alpha)\right|+A(| w_n|^\alpha)\right).
$$
Since by Lemma \ref{F1},
$$
A(| t|)\geq  c_0| t|^m, \,\,  \forall t\in [-1,1],
$$
it follows that
$$
\left(\int_{B_R(y)}| w_n|^{\frac{\alpha mN}{N-1}}\right)^{\frac{N}{N-1}}\leq C \int_{B_R(y)} \left(a(| w_n|)| w_n|| \nabla w_n|+A(| w_n|)\right).
$$
Next, let us fix $\alpha>0$ large enough and $p=\frac{m}{N}+m\alpha$. Thereby,
\begin{eqnarray*}
\int_{B_R(y)}| w_n|^p&=&\int_{B_R(y)}| w_n|^{\frac{m}{N}}|w_n|^{m\alpha}\\
&\leq& \left(\int_{B_R(y)}| w_n|^m\right)^{\frac{1}{N}}\left(\int_{B_R(y)}| w_n|^{\frac{m\alpha N}{N-1}}\right)^\frac{N-1}{N}.
\end{eqnarray*}
By (\ref{WN}),
$$
\left(\int_{B_R(y)}| w_n|^m\right)^{\frac{1}{N}}<\epsilon,
$$
for $n$ large enough and for all $y \in \mathbb{R}^{N}$. Hence, there is $n_0 \in \mathbb{N}$ such that
$$
\int_{B_R(y)}| w_n|^p\leq \epsilon c_1\int_{B_R(y)} f_n,\ n\geq n_0 \,\,\, \mbox{and} \,\,\,  y\in\mathbb{R}^N,
$$
where
$$
f_n=a(| w_n|)| w_n|| \nabla w_n|+A(| w_n|).
$$
Now, we set $\{y_j\}_{j\in\mathbb{N}}\subset\mathbb{R}^N$ such that $\mathbb{R}^N=\displaystyle\cup_{j\in\mathbb{N}} B_R(y_j)$ and each point of $\mathbb{R}^N$ is contained in at most $\kappa$ balls. Then,
\begin{eqnarray*}
\int_{\mathbb{R}^N} |w_n|^p&\leq& \displaystyle \sum_{j\in\mathbb{N}}\int_{B_R(y_j)}|w_n|^p\leq \epsilon c_1 \sum_{j\in\mathbb{N}}\int_{B_R(y_j)}f_n\\
&\leq& \displaystyle\epsilon c_1 \sum_{j\in\mathbb{N}}\int_{\mathbb{R^N}}f_n\chi_{B_R(y_j)}
\leq\displaystyle\epsilon c_1 \int_{\mathbb{R^N}}f_n\sum_{j\in\mathbb{N}}\chi_{B_R(y_j)}\\
&\leq& \epsilon c_1 \kappa \int_{\mathbb{R}^N}f_n.
\end{eqnarray*}
As $\{u_n \}$ is bounded in $W^{1}L_A(\mathbb{R}^{N})$, the sequence $\{f_n \}$ is bounded in $L^{1}(\mathbb{R}^{N})$. In this way, the last inequality gives
$$
w_n\buildrel  n \to +\infty \over
\longrightarrow0\ \text{in }L^p(\mathbb{R}^N),
$$
for  $p$ large enough. On the other hand,
$$
|w_n|^m_m=\int_{\mathbb{R^N}}| w_n|^m\leq c_0\int_{\mathbb{R}^N} A(| w_n|)\leq C, \ n\in\mathbb{N},
$$
from where it follows that $\{w_n\}$ is bounded in $L^m(\mathbb{R}^N)$. Then, by interpolation,
$$
w_n\buildrel  n \to +\infty \over
\longrightarrow0\ \text{in } L^q(\mathbb{R}^N), \forall \ q>m.
$$
From $(a_2)$, it follows that $l^{*}, m^* > m$, thus
\begin{equation} \label{wn}
w_n\buildrel  n \to +\infty \over
\longrightarrow0\ \text{in }L^{l^*}(\mathbb{R}^N)\ \text{and }L^{m^*}(\mathbb{R}^N).
\end{equation}
On the other hand, by Lemma \ref{F2},
$$
A_*(t)\leq C (|t|^{m^*}+|t|^{l^*}) \,\,\, \forall t \in \mathbb{R}^{N}.
$$
This combined with (\ref{wn}) gives
$$
\int_{\mathbb{R}^N}A_*(| w_n|) \to 0.
$$
From $(B_1)-(B_2)$, given $\epsilon >0$, there exists $C_{\epsilon}>0$ verifying
$$
B(|t|)\leq \epsilon A(| t|) +C_\epsilon A_*(| t|), \,\,\, t\in \mathbb{R}.
$$
Therefore,
\begin{eqnarray*}
\int_{\mathbb{R}^N}B(| w_n|)&\leq &\epsilon \int_{\mathbb{R}^N}A(| w_n|)+C_\epsilon\int_{\mathbb{R}^N}A_*(| w_n|)\\
&\leq& \epsilon C +C_\epsilon \int_{\mathbb{R}^N}A_*(| w_n|),
\end{eqnarray*}
from where it follows that
$$
\displaystyle \limsup_{n\rightarrow+\infty} \int_{\mathbb{R}^N}B(| w_n|)\leq \epsilon C \,\,\, \forall \epsilon>0,
$$
showing that
$$
\displaystyle \limsup_{n\rightarrow+\infty} \int_{\mathbb{R}^N}B(|w_n|)=0,
$$
that is,
$$
w_n\buildrel  n \to +\infty \over
\longrightarrow0\ \text{in }L_B(\mathbb{R}^N),
$$
finishing the proof of lemma.  \hspace{1 cm} $\rule{2mm}{2mm}$

\vspace{0.5 cm}

\noindent {\bf Proof of the compactness theorem  involving the group $O(N)$ ( see Theorem \ref{T5} )}

\vspace{0.5 cm}

The proof follows the same arguments used in Willem \cite[Theorem 1.24]{Willem}, when $|u_n|^{2}$ is replaced by $A(|u_n|)$. Here, we will make a sketch of the proof for convenience of the reader.

\vspace{0.5 cm}

Let $\{u_n\} $ be a sequence in  $W^{1}_{0,G}L_{A}(\Omega)$ with
$$
u_n \rightharpoonup 0 \,\,\, \mbox{in} \,\,\,  W^{1}_{0,G}L_{A}(\Omega).
$$
Without loss of generality, we can assume that $\{u_n \} \subset W^{1}_{0}L_{A}(\mathbb{R}^{N})$ by supposing that $u_n(x)=0$ for all $x \in \Omega^{c}$.

From definition of $m(y,r,G)$,
$$
\int_{B_{r}(y)}A(|u_n|) \leq  \frac{\displaystyle \sup_{n} \int_{\mathbb{R}^{N}}A(|u_n|)}{m(y,r,G)} \,\,\,\,\, \forall n \in \mathbb{N} \,\,\, \mbox{and} \,\,\, y \in \mathbb{R}^{N}.
$$
Once that $\Omega$ is compatible with $G$, given $\epsilon >0$, there is $R>0$ such that
\begin{equation} \label{W00}
\sup_{|y| \geq R}\int_{B_{r}(y)}A(|u_n|) \leq \epsilon, \,\,\, \forall n \in \mathbb{N}.
\end{equation}
On the other hand, one has
$$
B_r(y) \subset B_{R+r}(0) \,\,\ \forall y \in B_{R}(0)
$$
which implies that
\begin{equation} \label{W0}
\sup_{|y| < R}\int_{B_{r}(y)}A(|u_n|) \leq \int_{B_{R+r}(0)}A(|u_n|).
\end{equation}
By (\ref{M2}),
$$
u_n \to 0 \,\,\ \mbox{in} \,\,\ L_{A}(B_{R+r}(0))
$$
that is,
\begin{equation} \label{W1}
\int_{B_{R+r}(0)}A(|u_n|) \to 0.
\end{equation}
Thereby, from (\ref{W0}) and (\ref{W1}), there exists $n_0 \in \mathbb{N}$ such that
$$
\sup_{|y| < R}\int_{B_{r}(y)}A(|u_n|) \leq  \epsilon  \,\,\, \forall n \geq n_0.
$$
Hence, from (\ref{W00}) and (\ref{W1}),
$$
\displaystyle \lim_{n\rightarrow +\infty}\displaystyle \sup_{y\in{\mathbb{R}^N}}\int_{B_{R}(y)}A(|u_{n}|)=0.
$$
Now, the result follows applying the Theorem \ref{dolions}.  \hspace{1 cm} $\rule{2mm}{2mm}$

\section{Existence of solutions for problem (P)}

In this section, we will use the results obtained in the previous section to prove Theorem \ref{gio2}.  Hereafter, let us denote by
$J: X \to \mathbb{R}$ the energy functional related to $(P)$ given by
$$
J(u) =\displaystyle \int_{\mathbb{R}^N}A(|\nabla u|) +
          \displaystyle \int_{\mathbb{R}^N}V(x)A(|u|)
          -\displaystyle \int_{\mathbb{R}^N}F(u),
$$
where $X=W^{1}L_{A}(\mathbb{R}^{N})$ when $V$ is periodic and
$$
X=\left\{ u \in  W^{1}L_{A, rad}(\mathbb{R}^{N})\,; \, \int_{\mathbb{R}^{N}}V(x)A(|u|) < +\infty \right\}
$$
when $V$ is a radial function. In both cases, $X$ will be endowed with the norm
$$
\|u\|=\|\nabla u\|_{A} + \|u\|_{V,A}
$$
where
$$
\|u\|_{V,A}= \inf \biggl\{\alpha >0\,;\, \int_{\mathbb{R}^{N}}V(x)A\Big(\frac{|u|}{\alpha}\Big)\leq 1\biggl\}.
$$
A simple computation gives that the above norm is equivalent to the usual norm of $W^{1}L_{A}(\mathbb{R}^{N})$ when $V$ is a continuous periodic function satisfying $(V_1)$. Moreover, it is possible to prove that  $J \in C^{1}(X,\mathbb{R})$ with
$$
J'(u)\phi =\displaystyle \int_{\mathbb{R}^N}a(|\nabla u|)\nabla u \nabla
\phi +
          \displaystyle \int_{\mathbb{R}^N}V(x)a(|u|)u\phi
          -\displaystyle \int_{\mathbb{R}^N}f(u)\phi,
$$
for all $\phi\in X$. \\

Our goal is looking for critical points of $J$, because its critical points are weak solutions for $(P)$. Next, we will show three lemmas for the functional $J$, which are true when $V$ is radial or periodic. These lemmas will occur, because the below embeddings
$$
X \hookrightarrow L_{A}(\mathbb{R}^{N}) \,\,\ \mbox{and} \,\,\, X \hookrightarrow L_{A_*}(\mathbb{R}^{N})
$$
are continuous. The first of them establishes that $J$ verifies the mountain pass geometry on $X$.
\begin{lem}\label{passo} If $(a_1)-(a_2),(f_{1})-(f_{2})$ and $(V_{1})$ hold, the functional $J$ satisfies the following conditions:\\
(i) There exist $\rho$, $\eta>0$, such that $J(u)\geq
\eta$, if $\|u\|=\rho$. \\

\noindent (ii) For any $ \phi \in C^{\infty}_{0}(\mathbb{R}^{N}) \setminus \{0\}$, $J(t \phi)\rightarrow - \infty$ as $t\mapsto +\infty$.
\end{lem}
\textbf{Proof. (i)} From assumptions $(f_{1})-(f_{2})$, given $\epsilon >0$, there exists $C_\epsilon>0$ such that
\begin{eqnarray*}
0\leq f(t)t\leq \epsilon a(|t|)|t|^{2}+ C_{\epsilon}a_{*}(|t|)|t|^{2}\,\,\,\,\,\,\,\, \forall t \in \mathbb{R}.
\end{eqnarray*}
From $(a_2)$ and Lemma \ref{F3},
\begin{eqnarray}\label{crtescimentof}
0\leq f(t)t\leq \epsilon m A(|t|)+ C_\epsilon m^{*} A_{*}(|t|)\,\,\,\,\,\,\,\, \forall t \in \mathbb{R}.
\end{eqnarray}
Using $(f_{3})$,
\begin{eqnarray}\label{crtescimentoF}
0\leq F(t)\leq \displaystyle\frac{\epsilon m}{\theta} A(|t|)+ \widetilde{C} A_{*}(|t|)\,\,\,\,\,\,\,\, \forall t \in \mathbb{R}.
\end{eqnarray}
From (\ref{crtescimentoF}) and $(V_{1})$,
$$
J(u) \geq \displaystyle \int_{\mathbb{R}^N}A(|\nabla u|) +\left(1 - \frac{\epsilon m}{\theta V_0}\right)
          \displaystyle \int_{\mathbb{R}^N}V(x)A(|u|)
           - C\displaystyle \int_{\mathbb{R}^N}A_{*}(|u|).
$$
Hence, for $\epsilon$ small enough, the Lemmas \ref{F1} and \ref{F2} imply that
$$
J(u) \geq C_{1} \displaystyle \left(\xi_{0}(\|\nabla u\|_{A}) +\xi_{0}(\|u\|_{V,A})\right)
           - C_{2}\xi_{3}(\|u\|_{A_{*}}).
$$
Choosing $\rho >0$  such that
$$
\|u\|=\|\nabla u\|_{A} + \|u\|_{V,A} = \rho <1  \ \ \mbox{and} \ \ \|u\|_{A_{*}}\leq C(\|\nabla u \|_{A}+ \|u\|_{V,A}) <\rho <1,
$$
we obtain
$$
J(u) \geq C_{1} \displaystyle (\|\nabla u\|^{m}_{A} +\|u\|^{m}_{V,A})
           - C_{2}\|u\|^{l^{*}}_{A_{*}},
$$
which yields
$$
J(u) \geq C_{3}\|u\|^{m} - C_4 \|u\|^{l^{*}},
$$
for some positive constants $C_3$ and $C_4$.  Since $0<m<l^{*}$, there exists $\eta >0$ such that
$$
J(u)\geq \eta \,\,\, \mbox{ for all} \,\,\, \|u\|=\rho.
$$
\noindent \textbf{(ii)}  From $(f_{3})$, there exist $C_{5}, C_{6}>0$ such that
$$
F(t)\geq C_{5}|t|^{\theta} - C_{6}, \ \ \mbox{for all} \ \ t \in \mathbb{R}.
$$
Fixing $ \phi \in C^{\infty}_{0}(\mathbb{R}^{N}) \setminus \{0\}$, the last inequality leads to
$$
J(t \phi) \leq \xi_{1}(t) ( \xi_{1}(\| \nabla \phi \|_{A}) + \xi_{1}(\|  \phi \|_{V,A}))  - C_{5}t^{\theta}\displaystyle\int_{\mathbb{R}^{N}}|\phi|^{\theta} + C_{6} \mbox{supp} \phi.
$$
Thus, for $t$ sufficient large,
$$
J(t \phi) \leq t^{m}( \xi_{1}(\| \nabla \phi \|_{A}) + \xi_{1}(\|  \phi \|_{V,A}))- C_{5}t^{\theta}\displaystyle\int_{\mathbb{R}^{N}}|\phi|^{\theta} + C_{6} \mbox{supp} \phi.
$$
Since $ m < \theta$, the result follows. $\rule{2mm}{2mm}$

\vspace{0.5 cm}

Now, in view of the last lemma, we can apply a version of
Mountain Pass Theorem without the Palais-Smale condition found in
\cite{BN} to get a sequence $\{u_n\} \subset  X$ verifying
\begin{equation} \label {NPM}
J(u_{n})\rightarrow c
\,\,\, \mbox{and} \,\,\, J'(u_{n})\rightarrow 0 \,\,\, \mbox{as} \,\,\, n\rightarrow \infty,
\end{equation}
where the level $c$ is characterized by
$$
c =  \displaystyle \inf_{\gamma\in
\Gamma} \displaystyle \max_{t\in[0,1]}J(\gamma(t))>0
$$
and $\Gamma =
\{\gamma \in C([0,1],X): J(0)=0 \,\, \mbox{and} \,\,\, J(\gamma(1))<
0\}$.

\vspace{0.5 cm}

\begin{lem}\label{doisitens}
Let $\{v_{n}\}$ be a $(PS)_{d}$ sequence for $J$.  Then, $\{v_n\}$ is a bounded sequence in $X$.
\end{lem}
\textbf{Proof.} Since $\{v_n\}$ is a $(PS)_d$ sequence for the functional $J$, there is $C>0$ such that
\begin{eqnarray*}
C(1+\|v_{n}\|)\geq  J(v_{n})-\frac{1}{\theta}J'(v_{n})v_{n}, \,\,\, \forall n \in \mathbb{N}.
\end{eqnarray*}
From $(f_{3})$,
\begin{eqnarray*}
C(1+\|v_{n}\|)&\geq &
\left(\frac{\theta-m}{\theta}\right)\displaystyle\int_{\mathbb{R}^{N}}A(|\nabla v_{n}|)+V(x)A(|v_{n}|) \\
&\geq & \left(\frac{\theta-m}{\theta}\right)\biggl[\xi_{0}(\|\nabla v_{n}\|_{A})+ \xi_{0}(\|v_{n}\|_{V,A})\biggl].
\end{eqnarray*}

\noindent Suppose  for  contradiction that, up to a subsequence, $\|v_{n}\|\rightarrow +\infty$. This way, we need to study the following situations: \\

\noindent a) $\|\nabla v_{n}\|_{A}\rightarrow +\infty$ and $\|v_{n}\|_{V,A}\rightarrow +\infty$, \\

\noindent b) $\|\nabla v_{n}\|_{A}\rightarrow +\infty$ and $\|v_{n}\|_{V,A}$ is bounded, \\

\noindent and \\

\noindent c) $\|\nabla v_{n}\|_{A}$ is bounded and $\|v_{n}\|_{V,A}\rightarrow +\infty$. \\

In the first case, the Lemma \ref{F2} implies that
\begin{eqnarray*}
C(1+\|v_{n}\|)\geq  C_{1}\biggl[\|\nabla v_{n}\|^{l}_{A}+  \|v_{n}\|^{l}_{V,A}\biggl]\geq C_{2}\|v_{n}\|^{l},
\end{eqnarray*}
for $n$ large enough, which is an absurd.\\

In case b), we have for $n$ large enough
\begin{eqnarray*}
C_3(1+\|\nabla v_{n}\|_{A}) \geq C(1+\|v_{n}\|) \geq C_{2}\|\nabla v_{n}\|^{l}_{A},
\end{eqnarray*}
which is an absurd again. The last case is similar to the case b). $\rule{2mm}{2mm}$

\vspace{0.5 cm}

Using the fact that $X$ is reflexive, it follows from Lemma \ref{doisitens} that there exists a subsequence of $\{u_n\}$, still denoted by itself, and $u \in X$ such that
$$
u_n \rightharpoonup u \,\,\, \mbox{in} \,\,\,\, X.
$$

\begin{lem} \label{gradiente} The sequence $\{u_n \}$ satisfies the following limit
$$
\nabla u_n(x)\buildrel n \to +\infty \over \longrightarrow\nabla
u(x)\ \mbox{a.e in} \,\,\, \mathbb{R}^{N}.
$$
As a consequence, we deduce that $u$ is a critical point for $J$, that is, $J'(u)=0$.
\end{lem}
\noindent \textbf{Proof.} We begin this proof observing that $(a_1)$ yields
\begin{equation} \label{monot}
\left(a(| x| )x-a(| y|)y\right)(x-y)>0, \,\,\, \forall  x,y \in \mathbb{R}^{N} \,\,\, \mbox{with} \,\,\, x \not= y.
\end{equation}
Given $R>0$, let us consider $\xi=\xi_R\in
C_0^\infty(\mathbb{R}^N)$ satisfying
$$
0\leq\xi\leq 1, \xi \equiv 1 \,\,\, \mbox{in} \,\,\, B_R(0) \,\,\, \mbox{and} \,\,\, \text{supp}(\xi)\subset B_{2R}(0).
$$
Using the above information,
\begin{eqnarray}\label{lim0.1}
0&\leq& \int_{B_R(0)}\left(a(|\nabla u_n|)\nabla u_n-a(|
\nabla u|)\nabla u\right)(\nabla u_n-\nabla u)\nonumber\\
&&\leq\int_{B_{2R}(0)}\left( a(| \nabla u_n|)\nabla
u_n-a(|\nabla u|)\nabla u\right)(\nabla u_n-\nabla u)\xi\\
&&=\int_{B_{2R}(0)}a(|\nabla u_n|)\nabla u_n(\nabla u_n-\nabla
u)\xi-\int_{B_{2R}(0)}a(|\nabla u|)\nabla u(\nabla u_n-\nabla
u)\xi.\nonumber
\end{eqnarray}
Now, combining the boundedness of $\{(u_n-u)\xi\}$  in $X$ with the limit $\|J'(u_n)\|=o_n(1)$, it follows that
\begin{eqnarray}\label{lim0.2}
o_n(1)&=&\int_{B_{2R}(0)} a(| \nabla u_n|)\nabla
u_n\nabla\left((u_n-u)\xi\right)+\int_{B_{2R}(0)}V(x)a(|u_n|)u_n(u_n-u)\xi\nonumber\\
&&-\int_{B_{2R}(0)}f(u_n)(u_n-u)\xi.
\end{eqnarray}
Note that $\{a(| u_n|)u_n\}$ is bounded in $L_{\widetilde{A}}(B_{2R}(0))$, because
$$
\int_{B_{2R}(0)}\widetilde{A}\left(a(| u_n |)u_n\right) \leq
\int_{B_{2R}(0)}A(2| u_n|)\leq K \int_{B_{2R}(0)}A(|
u_n|)<+\infty.
$$
From this
\begin{eqnarray}\label{lim0.3}
\left|\int_{B_{2R}(0)}V(x)a(| u_n|)u_n(u_n-u)\xi\right|&\leq&
\int_{B_{2R}(0)}|V(x)|| a(| u_n|) u_n| |u_n-u|\nonumber\\
&\leq&2M \parallel a(| u_n|)|
u_n|\parallel_{\widetilde{A},B_{2R}(0)}\parallel
u_n-u\parallel_{A,B_{2R}(0)}\nonumber\\
&\leq& C_1 \parallel u_n-u\parallel_{A,B_{2R}(0)}\to 0.
\end{eqnarray}
where $M=\displaystyle \sup_{x \in B_{2R}(0)}|V(x)|$. On the other hand, using again the boundedness of $(u_n)$ in $X$ and (\ref{D1}),
$$
\int_{B_{2R}(0)}\widetilde{A}_*(a_*(u_n)u_n)\leq
\int_{B_{2R}(0)}A_*(2u_n)\leq C_2, \ n\in \mathbb{N},
$$
implying that $\{a_*(u_n)u_n\}$ is bounded in $L_{\widetilde{A}_*}(B_{2R}(0))$.  Since
\begin{eqnarray}\label{lim1}
\left|\int_{B_{2R}(0)}f(u_n)(u_n-u)\xi\right|&\leq&
 \epsilon \Big(\int_{B_{2R}(0)}|
a(|
u_n|)u_n|| u_n-u|\nonumber\\
&&+\int_{B_{2R}(0)}| a_*(u_n)u_n|| u_n-u|\Big) +
c_1\int_{B_{2R}(0)} |u_n-u|\nonumber\\
&\leq& \epsilon\parallel a(|
u_n|)u_n\parallel_{\widetilde{A},B_{2R}(0)}\parallel
u_n-u\parallel_{A,B_{2R}(0)}\\
&&+\epsilon c_2\parallel a_*(u_n)u_n
\parallel_{\widetilde{A}_*,B_{2R}(0)}\parallel
u_n-u\parallel_{A_*,B_{2R}(0)}\nonumber\\
&&+c_3\parallel u_n-u\parallel_{A,B_{2R}(0)}\nonumber,
\end{eqnarray}
the boundedness of $\{u_n\}, \{a(u_n)u_n\}$ and
$\{a_*(u_n)u_n\}$ in  $L_A(B_{2R}(0)),$ 
$L_{\widetilde{A}}(B_{2R}(0))$ and  $L_{\widetilde{A}_*}(B_{2R}(0))$ respectively lead to
$$
\left|\int_{B_{2R}(0)}f(u_n)(u_n-u)\xi\right| \leq \epsilon C_4 +c_3\parallel u_n-u\parallel_{A,B_{2R}(0)}.
$$
Now, using the convergence of  $\{u_n\}$ to $u$ in $L_A(B_{2R}(0))$, we get
\begin{equation}\label{lim2}
\left|\int_{B_{2R}(0)}f(u_n)(u_n-u)\xi\right|\to 0.
\end{equation}
A similar idea can be used to establish the limit
\begin{eqnarray}
\int_{B_{2R}(0)} (u_n-u) a(|\nabla u_n|)\nabla
u_n\nabla \xi \to 0.
\end{eqnarray}
Moreover, the weak convergence of $\{u_n\}$ to $u$ in $W^1L_{A}(\mathbb{R}^N)$ gives
\begin{equation}\label{lim3}
\int_{B_{2R}(0)}\xi a(|\nabla u|)\nabla u(\nabla u_n-\nabla
u)\to 0.
\end{equation}
From (\ref{lim0.1})-(\ref{lim3}),
$$
\int_{B_R(0)}\left(a(|\nabla u_n|)\nabla u_n-a(|
\nabla u|)\nabla u\right)(\nabla u_n-\nabla u)\to0.
$$
Setting $\beta:{\mathbb{R}}^N\rightarrow{\mathbb{R}} ^N$ by
$$
\beta(x)=a(| x|) x, \ x\in{\mathbb{R}}^N,
$$
and observing that $\beta$ is monotone by (\ref{monot}), the last limit imply that for some subsequence, still denoted by itself,
$$
\left(\beta(\nabla u_n(x))-\beta(\nabla u(x))\right)(\nabla
u_n(x)-\nabla u(x))\to 0 \ \mbox{a.e in  } \,\,\, B_{R}(0).
$$
Applying a result found in Dal Maso and Murat \cite{Maso}, it follows that
$$
\nabla u_n(x)\to \nabla u(x) \ \mbox{a.e in } \,\,\, B_{R}(0),
$$
for each $R>0.$ As $R$ is arbitrary,  there is a subsequence of $\{u_n\}$, still denoted by itself, such that
$$
\nabla u_n(x)\to \nabla u(x) \ \mbox{a.e in} \,\,\, \mathbb{R}^N.
$$
Recalling that $\{a(|\nabla u_n|)\frac{\partial u_n}{\partial x_i}\}$ is bounded in $L_{\tilde{A}}(\mathbb{R}^{N})$, we get from Lemma \ref{brezislieb}
$$
\int_{\mathbb{R}^{N}}a(|\nabla u_n|)\nabla u_n \nabla v \to \int_{\mathbb{R}^{N}}a(|\nabla u|)\nabla u \nabla v,
$$
for all $v \in X_{c}=\{v \in X\,\;\, \,\,  v \,\, \mbox{has compact support} \,\, \}$. On the other hand, once that $V$ is bounded on the support of $v$, $\{a(| u_n|)u_n \}$ is bounded in $L_{\tilde{A}}(\mathbb{R}^{N})$ and $\{a_*(| u_n|)u_n \}$ is bounded in $L_{\tilde{A}_*}(\mathbb{R}^{N})$, we have again by Lemma  \ref{brezislieb}
$$
\int_{\mathbb{R}^{N}}V(x)a(| u_n|)u_n v \to \int_{\mathbb{R}^{N}}V(x)a(| u|)u v
$$
and
$$
\int_{\mathbb{R}^{N}}f(u_n)v \to \int_{\mathbb{R}^{N}}f(u)v.
$$
Therefore,
$$
J'(u)v=0 \,\,\, \forall v \in X_{c}.
$$
Now, the lemma follows using the fact that $X_c$ is dense in $X$. \hspace{1 cm} \rule{2mm}{2mm}

\vspace{0.5 cm}

\noindent \subsection{Proof of Theorem \ref{gio2}}

The reader is invited to observe that the main difference between the radial and periodic case is the following: In the radial case, the Theorem \ref{compacidade} permits to prove that the energy functional $J$ verifies the $(PS)$ condition, while in the periodic case, we do not have this condition and we overcome this difficulty by using the Theorem \ref{dolions}.

We will prove the Theorem \ref{gio2} studying firstly the radial case, and after, the periodic case.

\vspace{0.5 cm}

\noindent {\bf The radial case:}

\vspace{0.5 cm}

For the radial case, we begin showing the following lemma

\begin{claim} \label{Convergencia de f}
Let $\{u_n\}$ the sequence given in (\ref{NPM}). If $(f_1)-(f_2)$ hold, one have
$$
\int_{\mathbb{R}^{N}}f(u_n)u_n \to \int_{\mathbb{R}^{N}}f(u)u.
$$
\end{claim}
Indeed, as  $\{u_n\}$ is a bounded sequence in $W^{1}L_{A, rad}(\mathbb{R}^{N})$,
$$
\sup_{n \in  \mathbb{N}}\int_{\mathbb{R}^{N}}(A_*(u_n) + A(u_n)) < +\infty.
$$
Moreover, by hypotheses  $(f_1)-(f_2)$, the function $P(t)=f(t)t$ verifies the limit
$$
\lim_{|t| \to 0} \frac{P(t)}{A(t)+A_{*}(t)}=0 \,\,\, \mbox{and} \,\,\, \lim_{|t| \to +\infty} \frac{P(t)}{A(t)+A_{*}(t)}=0.
$$
Since by Theorem \ref{strauss},
$$
u_n(x) \to 0 \,\,\, \mbox{as} \,\,\, |x| \to +\infty, \,\,\, \mbox{uniformly with respect to} \,\, n,
$$
it follows from \cite[Theorem A.I]{bl},
$$
\int_{\mathbb{R}^{N}}f(u_n)u_n \to \int_{\mathbb{R}^{N}} f(u)u,
$$
proving the claim.

\vspace{0.5 cm}

Recalling that $J'(u_n)u_n =o_n(1)$, or equivalently,
$$
\int_{\mathbb{R}^{N}}(a(|\nabla u_n|)|\nabla u_n|^{2} + V(x)a(|u_n|)|u_n|^{2})=\int_{\mathbb{R}^{N}}f(u_n)u_n +o_n(1),
$$
we derive from Claim \ref{Convergencia de f}
$$
\lim_{n \to \infty}\int_{\mathbb{R}^{N}}(a(|\nabla u_n|)|\nabla u_n|^{2} + V(x)a(|u_n|)|u_n|^{2})=\int_{\mathbb{R}^{N}}f(u)u.
$$
Using the fact that $J'(u)u=0$, it follows that
$$
\lim_{n \to \infty}\int_{\mathbb{R}^{N}}(a(|\nabla u_n|)|\nabla u_n|^{2} + V(x)a(|u_n|)|u_n|^{2})=\int_{\mathbb{R}^{N}}(a(|\nabla u|)|\nabla u|^{2} + V(x)a(|u|)|u|^{2}).
$$
Once that
$$
\nabla u_n(x) \to \nabla u(x) \;\;\; \mbox{and} \;\;\; u_n(x) \to u(x) \;\;\; \mbox{a.e in} \;\;\; \mathbb{R}^{N},
$$
we conclude that
$$
a(|\nabla u_n|)|\nabla u_n|^{2} \to a(|\nabla u|)|\nabla u|^{2} \;\;\; \mbox{in} \;\;\; L^{1}(\mathbb{R}^{N})
$$
and
$$
V(x)a(|u_n|)|u_n|^{2} \to V(x)a(|u|)|u|^{2} \;\;\; \mbox{in} \;\;\; L^{1}(\mathbb{R}^{N}).
$$
These limits combined with $(a_2)$ yields
$$
\int_{\mathbb{R}^{N}}A(|\nabla u_n - \nabla u|) \to 0
$$
and
$$
\int_{\mathbb{R}^{N}}V(x)A(|u_n - u|) \to 0.
$$
Hence, by a similar arguments used in (\ref{CV1}), we derive that
$$
u_n \to u \;\;\; \mbox{in} \;\;\; X,
$$
and thus,
$$
J(u)=c>0 \;\;\; \mbox{and} \;\;\; J'(u)=0,
$$
showing that $u$ is a critical point of $J$ in  $X$. Now, using a principle of symmetric criticality on reflexive Banach spaces due to \linebreak de Morais Filho, Do \'O and Souto \cite{DMO}, we have that  $u$ is a critical point of $J$ in  $W^{1}L_{A}(\mathbb{R}^N)$, and so, $u$ is a  nontrivial solution for problem $(P)$.

\vspace{0.5 cm}

\noindent {\bf The periodic case}

\vspace{0.5 cm}

By Lemma \ref{gradiente}, we know that the weak limit $u$ of the sequence $\{u_n\}$ given in (\ref{NPM}) is a critical point for $J$. If $u\not= 0$, the theorem is proved. However, if $u = 0$, we have the following claim:\\

\begin{claim} \label{afirmacao2}
\, There is $R> 0$ such that
\begin{eqnarray}\label{udiferentedezero}
\displaystyle \liminf_{n\rightarrow
+\infty}\displaystyle \sup_{y\in{\mathbb{R}^N}}\int_{B_{R}(y)}A(u_{n}) >0.
\end{eqnarray}
\end{claim}
In fact, if the above claim does not hold, by using
Theorem \ref{dolions}, we derive the limit
\begin{equation} \label{LW}
\displaystyle \int_{{\mathbb{R}^N}}B(|u_{n}|)\rightarrow 0,
\end{equation}
for any $N$-function $B$ satisfying $(B_1)-(B_2)$. Fixing a $N$-function $B$ satisfying  $(B_1)-(B_2)$, it follows from $(f_1)-(f_2)$ that given $\epsilon >0$, there exists $C_\epsilon >0$ such that
$$
|f(u_n)u_n| \leq \epsilon ( A(|u_n|) +A_{*}(|u_n|) ) + C_\epsilon B(|u_n|) \,\,\, \forall n \in \mathbb{N}.
$$
Thereby, the above inequality together with (\ref{LW}) gives
$$
\int_{\mathbb{R}^{N}}f(u_n)u_n \to 0.
$$
Recalling that $J'(u_n)u_n=o_n(1)$, that is,
$$
\int_{\mathbb{R}^{N}}a(|\nabla u_n|)|\nabla u_n|^{2}+\int_{\mathbb{R}^{N}}V(x)a(|u_n|)|u_n|^{2}=\int_{\mathbb{R}^{N}}f(u_n)u_n + o_{n}(1),
$$
we obtain
$$
\int_{\mathbb{R}^{N}}a(|\nabla u_n|)|\nabla u_n|^{2}+\int_{\mathbb{R}^{N}}V(x)a(|u_n|)|u_n|^{2} \to 0.
$$
The last limit together with $(a_2)$ gives
$$
\int_{\mathbb{R}^{N}}A(|\nabla u_n|)+\int_{\mathbb{R}^{N}}V(x)A(|u_n|)  \to 0,
$$
implying that $\{u_{n}\}$ converges strongly to zero in $W^1L_{A}(\mathbb{R}^N)$, leading to  $c = 0$, which is an absurd. Thus, the limit (\ref{udiferentedezero}) holds and the claim is proved.

\vspace{0.5 cm}

Therefore, there are $R,\alpha>0$ and $\{y_{n}\} \subset \mathbb{Z}^{N}$ such that
\begin{eqnarray}\label{udiferentedezero1}
\int_{B_{R}(y_n)}A(u_{n})>\alpha.
\end{eqnarray}
Now, letting
$\overline{u}_{n}(x) = u_{n}(x-y_{n})$, since $V$ is $\mathbb{Z}^{N}$-periodic
function, one has
$$
\|\overline{u}_{n}\| =\|u_{n}\|, \,\, J({\overline u}_{n}) = J(u_{n}) \,\,\, \mbox{and} \,\,\,
J'({\overline u}_{n})=o_n(1).
$$
Then, there exists ${\overline u} $ such that ${\overline u_{n}}\rightharpoonup {\overline
u}$ weakly in $W^1L_{A}(\mathbb{R}^N)$, and as before, it follows that
$J'({\overline u})= 0$. Now, by (\ref{udiferentedezero1}),
$$
\int_{B_{R}(0)}A({\overline u}_{n})\geq  \alpha >0,
$$
which together with the compact embeddings yields
$$
\int_{B_{R}(0)}A({\overline u})\geq  \alpha >0,
$$
showing that ${\overline u} \not= 0$, and thereby, finishing the proof of the Theorem \ref{gio2}. $\rule{2mm}{2mm}$

\vspace{1 cm}

\noindent \textbf{Acknowledgments.} The authors are  grateful to the referees
for a number of helpful comments for improvement in this article.

\end{document}